\documentclass[11pt]{article}
\usepackage{amssymb}
\usepackage{amsfonts}
\usepackage{amsmath}

\input{tcilatex}
\begin{document}

\begin{center}
\textbf{A numerical method based on the reproducing kernel Hilbert space
method for the solution of fifth-order boundary-value problems}

Mustafa Inc, Ali Akg\"{u}l and Mehdi Dehghan

Department of Mathematics, Science Faculty, F\i rat University, 23119 Elaz\i 
\u{g} / T\"{u}rkiye

Department of Mathematics, Education Faculty, Dicle University, 21280
Diyarbak\i r / T\"{u}rkiye

Department of Applied Mathematics, Faculty of Mathematics and Computer
Science, Amirkabir University of Technology, No. 424, Hafez Ave., Tehran,
Iran

minc@firat.edu.tr
\end{center}

\textbf{Abstract:} In this paper, we present a fast and accurate numerical
scheme for the solution of fifth-order boundary-value problems. We apply the
reproducing kernel Hilbert space method (RKHSM) for solving this problem.
The analytical results of the equations have been obtained in terms of
convergent series with easily computable components. We compare our results
with spline methods, decomposition method, variational iteration method,
Sinc-Galerkin method and homotopy perturbation methods. The comparison of
the results with exact ones is made to confirm the validity and efficiency.


\bigskip

\textbf{1. Introduction}

In this work we consider the numerical approximation for the fifth-order
boundary-value problems of the form%
\begin{equation}
y^{\left( v\right) }=f\left( x\right) y+g\left( x\right) ,\quad x\in \left[
a,b\right] ,  \tag{1.1}
\end{equation}%
with boundary conditions%
\begin{equation}
y\left( a\right) =A_{0},\quad y^{\prime }\left( a\right) =A_{1},\quad
y^{\prime \prime }\left( a\right) =A_{2},\quad y\left( b\right) =B_{0},\quad
y^{\prime }\left( b\right) =B_{1},\quad  \tag{1.2}
\end{equation}%
where the functions $f\left( x\right) $ and $g\left( x\right) $ are
continuous on $\left[ a,b\right] $ and $A_{0},$ $A_{1},$ $A_{2},$ $B_{0},$ $%
B_{1}$ are finite real constants. For more details about computational code
of boundary value problems, the reader is referred to [1-3].

This type of boundary-value problems arise in the mathematical modelling of
viscoelastic flows and other branches of mathematical, physical and
engineering sciences [4,5] and references therein. Theorems which list the
conditions for the existence and uniqueness of solutions of such problems
are thoroughly discussed in a book by Agarwal [6]. Khan [7] investigated the
fifth-order boundary-value problems by using finite difference methods.
Wazwaz [8] applied Adomian decomposition method for solution of such type of
boundary-value problems. The use of spline function in the context of
fifth-order boundary-value problems was studied by Fyfe [9], who used the
quintic polynomial spline functions to develop consistency relations
connecting the values of solution with fifth-order derivatives at the
respective nodes. Polynomial sextic spline functions were used [10] to
develop the smooth approximations to the solution of problems (1.1) and
(1.2). Caglar et al. [11] have used sixth-degree B-spline functions to
develop first-order accurate method for the solution two-point special
fifth-order boundary-value problems. Noor and Mohyud-Din [12,13] applied
variational iteration and homotopy perturbation methods for solving the
problems (1.1) and (1.2), respectively. Khan [14] have used the
non-polynomial sextic spline functions for the solution fifth-order
boundary-value problems. El-Gamel [15] employed the sinc-Galerkin method to
solve the problems (1.1) and (1.2). Lamnii et al. [16] developed and
analyzed two sextic spline collocation methods for the problem. Siddiqi et
al. [17,18] used the non-polynomial sextic spline method for special
fifth-order problems (1.1) and (1.2). Wang et al. [19] attempted to obtain
upper and lower approximate solutions of such problems by applying the
sixth-degree B-spline residual correction method.

In this paper, the RKHSM [20,21] will be used to investigate the fifth-order
boundary-value problems. In recent years, a lot of attetion has been devoted
to the study of RKHSM to investigate various scientific models. The RKHSM
which accurately computes the series solution is of great interest to
applied sciences. The method provides the solution in a rapidly convergent
series with components that can be elegantly computed. The efficiency of the
method was used by many authors to investigate several scientific
applications. Geng and Cui [22] applied the RKHSM to handle the second-order
boundary value problems. Yao and Cui [23] and Wang et al. [24] investigated
a class of singular boundary value problems by this method and the obtained
results were good. Zhou et al. [25] used the RKHSM effectively to solve
second-order boundary value problems. In [26], the method was used to solve
nonlinear infinite-delay-differential equations. Wang and Chao [27], Li and
Cui [28], Zhou and Cui [29] independently employed the RKHSM to
variable-coefficient partial differential equations. Geng and Cui [30], Du
and Cui [31] investigated the approximate solution of the forced Duffing
equation with integral boundary conditions by combining the homotopy
perturbation method and the RKHSM. Lv and Cui [32] presented a new algorithm
to solve linear fifth-order boundary value problems. In [33,34], authors
developed a new existence proof of solutions for nonlinear boundary value
problems. Cui and Du [35] obtained the representation of the exact solution
for the nonlinear Volterra-Fredholm integral equations by using the
reproducing kernel Hilbert space method. Wu and Li [36] applied iterative
reproducing kernel method to obtain the analytical approximate solution of a
nonlinear oscillator with discontinuties. Recently, the method was apllied
the fractional partial differential equations and multi-point boundary value
problems [34-37]. For more details about RKHSM and the modified forms and
its effectiveness, see [20-43] and the references therein.

The paper is organized as follows. Section 2 is devoted to several
reproducing kernel spaces and a linear operator is introduced. S{\normalsize %
olution represantation in}\textbf{\ }$W_{2}^{6}[a,b]$ has been presented in
Section 3. It provides the main results, the exact and approximate solution
of $(1.1)$ and an iterative method are developed for the kind of problems in
the reproducing kernel space. We have proved that the approximate solution
converges to the exact solution uniformly. Some numerical experiments are
illustrated in Section 4. There are some conclusions in the last section.

\bigskip

\textbf{2. Preliminaries}

\textbf{2.1. Reproducing Kernel Spaces}\ 

In this section, we define some useful reproducing kernel spaces.

\textbf{Definition 2.1.} \textit{(Reproducing kernel)}. Let $E$ be a
nonempty abstract set. A function $K:E\times E\longrightarrow C$ is a
reproducing kernel of the Hilbert space $H$ if and only if%
\begin{equation}
\left\{ 
\begin{array}{c}
\forall t\in E,\text{ }K\left( .,t\right) \in H, \\ 
\forall t\in E,\text{ }\forall \varphi \in H,\text{ }\left\langle \varphi
\left( .\right) ,K\left( .,t\right) \right\rangle =\varphi \left( t\right) .%
\end{array}%
\right.  \tag{2.1}
\end{equation}

The last condition is called "the reproducing property": the value of the
function $\varphi $ at the point $t$ is reproduced by the inner product of $%
\varphi $ with $K\left( .,t\right) $

\bigskip

\textbf{Definition 2.2.}%
\begin{equation*}
W_{2}^{6}[0,1]=\left\{ 
\begin{array}{c}
u(x)\mid u(x),\text{ }u^{\prime }(x),\text{ }u^{\prime \prime }(x),\text{ }%
u^{\prime \prime \prime }(x),\text{ }u^{(4)}(x),\text{ }u^{(5)}(x)\text{\ }
\\ 
\text{are absolutely continuous in }[0,1], \\ 
u^{\left( 6\right) }(x)\in L^{2}[0,1],\text{ }x\in \lbrack 0,1],\text{ }%
u(0)=u(1)=u^{\prime }(0)=u^{\prime }(1)=0=u^{\prime \prime }(0)=0%
\end{array}%
\right\} ,
\end{equation*}%
The sixth derivative of $u(x)$ exists almost everywhere since $u^{(5)}(x)$
is absolutely continuous. The inner product and the norm in $W_{2}^{6}[0,1]$
are defined respectively by%
\begin{equation*}
\left\langle u(x),g(x)\right\rangle
_{W_{2}^{6}}=\sum_{i=0}^{5}u^{(i)}(0)g^{(i)}(0)+%
\int_{0}^{1}u^{(6)}(x)g^{(6)}(x)dx,\text{ \ }u(x),g(x)\in W_{2}^{6}[0,1],
\end{equation*}%
and

\begin{equation*}
\left\Vert u\right\Vert _{W_{2}^{6}}=\sqrt{\left\langle u,u\right\rangle
_{_{W_{2}^{6}}}},\ u\in W_{2}^{6}[0,1].
\end{equation*}%
The space $W_{2}^{6}[0,1]$ \ is a reproducing kernel space, i.e., for each
fixed $y\in \lbrack 0,1]$ \ and any $u(x)\in W_{2}^{6}[0,1],$ there exists a
function $R_{y}(x)$ such that

\begin{equation*}
u(y)=\left\langle u(x),\ R_{y}(x)\right\rangle _{W_{2}^{6}}.
\end{equation*}

\bigskip

\textbf{Definition 2.3.}\ 

\begin{equation*}
W_{2}^{1}[0,1]=\left\{ 
\begin{array}{c}
u(x)\mid u(x)\text{ is absolutely continuous in \ }[0,1]\text{ } \\ 
u^{\prime }(x)\in L^{2}[0,1],\text{ }x\in \lbrack 0,1],%
\end{array}%
\right\} ,
\end{equation*}%
The inner product and the norm in $W_{2}^{1}[0,1]$ \ are defined
respectively by

\begin{equation*}
\left\langle u(x),g(x)\right\rangle
_{W_{2}^{1}}=u(0)g(0)+\int_{0}^{1}u^{\prime }(x)g^{\prime }(x)dx,\text{ \ }%
u(x),g(x)\in W_{2}^{1}[0,1],
\end{equation*}%
and

\begin{equation*}
\left\Vert u\right\Vert _{W_{2}^{1}}=\sqrt{\left\langle u,u\right\rangle _{_{%
{\LARGE W}_{2}^{1}}}},\text{ }u\in W_{2}^{1}[0,1].
\end{equation*}%
The space $W_{2}^{1}[0,1]$ \ is a reproducing kernel space and its
reproducing kernel function ${\normalsize T}_{x}{\normalsize (y)}$ is given
by

\begin{equation}
{\normalsize T}_{x}{\normalsize (y)}{\LARGE =}\left\{ 
\begin{array}{c}
1+x,\text{ \ \ }x\leq y, \\ 
\\ 
1+y,\text{ \ \ }x>y.%
\end{array}%
\right.  \tag{2.2}
\end{equation}

\bigskip

\textbf{Theorem 2.1. }The space $W_{2}^{6}[0,1]$ is a reproducing kernel
Hilbert space whose reproducing kernel function is given by,%
\begin{equation*}
R_{y}\left( x\right) =\left\{ 
\begin{array}{c}
\sum_{i=1}^{12}c_{i}\left( y\right) x^{i-1},\quad x\leq y, \\ 
\\ 
\sum_{i=1}^{12}d_{i}\left( y\right) x^{i-1},\quad x>y,%
\end{array}%
\right.
\end{equation*}%
where

\begin{equation*}
c_{1}(y)=0,
\end{equation*}

\begin{equation*}
{\normalsize c}_{2}{\normalsize (y)}{\normalsize =0,}
\end{equation*}

\begin{equation*}
c_{3}(y)=0,
\end{equation*}

\begin{eqnarray*}
c_{4}(y) &=&\frac{2461}{42301989}y^{5}+\frac{2461}{253811934}y^{6}+\frac{335%
}{507623868}y^{7} \\
&& \\
&&-\frac{8321}{10660101228}y^{10}-\frac{7325}{1776683538}y^{8}-\frac{11725}{%
84603978}y^{4} \\
&& \\
&&+\frac{12221}{169207956}y^{3}+\frac{56255}{21320202456}y^{9}+\frac{2003}{%
21320202456}y^{11},
\end{eqnarray*}

\begin{eqnarray*}
c_{5}(y) &=&-\frac{158419}{1353663648}y^{5}-\frac{158419}{8121981888}y^{6}+%
\frac{99305}{14213468304}y^{7} \\
&& \\
&&+\frac{157481}{341123239296}y^{10}-\frac{335}{16243963776}y^{8}+\frac{%
728021}{2707327296}y^{4} \\
&& \\
&&-\frac{11725}{84603978}y^{3}-\frac{196265}{170561619648}y^{9}-\frac{2687}{%
42640404912}y^{11},
\end{eqnarray*}

\begin{eqnarray*}
c_{6}(y) &=&\frac{4056701}{67683182400}y^{5}-\frac{2011}{126905967}y^{6}+%
\frac{158419}{284269366080}y^{7} \\
&& \\
&&+\frac{11843}{304574320800}y^{10}+\frac{2461}{284269366080}y^{8}-\frac{%
158419}{1353663648}y^{4} \\
&& \\
&&+\frac{2461}{42301989}y^{3}-\frac{168263}{1705616196480}y^{9}-\frac{8999}{%
1705616196480}y^{11},
\end{eqnarray*}

\begin{eqnarray*}
c_{7}(y) &=&\frac{4056701}{406099094400}y^{5}-\frac{2011}{7614358020}y^{6}+%
\frac{158419}{1705616196480}y^{7} \\
&& \\
&&+\frac{11843}{1827445924800}y^{10}+\frac{2461}{1705616196480}y^{8}-\frac{%
158419}{8121981888}y^{4} \\
&& \\
&&+\frac{2461}{253811934}y^{3}-\frac{168263}{10233697178880}y^{9}-\frac{8999%
}{10233697178880}y^{11},
\end{eqnarray*}

\begin{eqnarray*}
c_{8}(y) &=&\frac{158419}{284269366080}y^{5}+\frac{158419}{1705616196480}%
y^{6}-\frac{19861}{596965668768}y^{7} \\
&& \\
&&-\frac{157481}{71635880252160}y^{10}+\frac{67}{682246478592}y^{8}-\frac{%
104003}{81219818880}y^{4} \\
&& \\
&&+\frac{335}{507623868}y^{3}+\frac{39253}{7163588025216}y^{9}+\frac{2687}{%
8954485031520}y^{11},
\end{eqnarray*}

\begin{eqnarray*}
c_{9}(y) &=&\frac{2461}{284269366080}y^{5}+\frac{2461}{1705616196480}y^{6}+%
\frac{67}{682246478592}y^{7} \\
&& \\
&&-\frac{8321}{71635880252160}y^{10}-\frac{1465}{2387862675072}y^{8}-\frac{%
335}{16243963776}y^{4} \\
&& \\
&&+\frac{12221}{1137077464320}y^{3}+\frac{11251}{28654352100864}y^{9}+\frac{%
2003}{143271760504320}y^{11},
\end{eqnarray*}

\begin{eqnarray*}
c_{10}(y) &=&-\frac{168263}{1705616196480}y^{5}-\frac{168263}{10233697178880}%
y^{6}+\frac{39253}{7163588025216}y^{7} \\
&& \\
&&+\frac{38153}{85963056302592}y^{10}+\frac{11251}{28654352100864}y^{8}-%
\frac{196265}{170561619648}y^{4} \\
&& \\
&&+\frac{56255}{21320202456}y^{3}-\frac{6313}{5372691018912}y^{9}-\frac{12751%
}{214907640756480}y^{11}-\frac{1}{725760}y^{2},
\end{eqnarray*}

\begin{eqnarray*}
c_{11}(y) &=&\frac{11843}{304574320800}y^{5}+\frac{11843}{1827445924800}%
y^{6}-\frac{157481}{71635880252160}y^{7} \\
&& \\
&&-\frac{45611}{268634550945600}y^{10}-\frac{8321}{71635880252160}y^{8}+%
\frac{157481}{341123239296}y^{4} \\
&& \\
&&-\frac{8321}{10660101228}y^{3}+\frac{38153}{85963056302592}y^{9}+\frac{%
49001}{2149076407564800}y^{11}+\frac{1}{3628800}y,
\end{eqnarray*}

\begin{eqnarray*}
c_{12}(y) &=&-\frac{8999}{1705616196480}y^{5}-\frac{8999}{10233697178880}%
y^{6}+\frac{2687}{8954485031520}y^{7} \\
&& \\
&&+\frac{49001}{2149076407564800}y^{10}+\frac{2003}{143271760504320}y^{8}-%
\frac{2687}{42640404912}y^{4} \\
&& \\
&&+\frac{2003}{21320202456}y^{3}-\frac{12751}{214907640756480}y^{9}-\frac{725%
}{236398404832128}y^{11}-\frac{1}{39916800},
\end{eqnarray*}

\begin{equation*}
d_{1}(y)=-\frac{1}{39916800}y^{11},
\end{equation*}

\begin{equation*}
d_{2}(y)=\frac{1}{3628800}y^{10},
\end{equation*}

\begin{equation*}
d_{3}(y)=-\frac{1}{725760}y^{9},
\end{equation*}

\begin{eqnarray*}
\text{\ }d_{4}(y) &=&\frac{12221}{169207956}y^{3}+\frac{12221}{1137077464320}%
y^{8}+\frac{2461}{42301989}y^{5} \\
&& \\
&&+\frac{2461}{253811934}y^{6}+\frac{335}{507623868}y^{7}-\frac{8321}{%
10660101228}y^{10} \\
&& \\
&&-\frac{11725}{84603978}y^{4}+\frac{56255}{21320202456}y^{9}+\frac{2003}{%
21320202456}y^{11},
\end{eqnarray*}

\begin{eqnarray*}
d_{5}(y) &=&\frac{728021}{2707327296}y^{4}-\frac{104003}{81219818880}y^{7}-%
\frac{158419}{1353663648}y^{5} \\
&& \\
&&-\frac{158419}{8121981888}y^{6}+\frac{157481}{341123239296}y^{10}-\frac{335%
}{16243963776}y^{8} \\
&& \\
&&-\frac{11725}{84603978}y^{3}-\frac{196265}{170561619648}y^{9}-\frac{2687}{%
42640404912}y^{11},
\end{eqnarray*}

\begin{eqnarray*}
d_{6}(y) &=&\frac{4056701}{67683182400}y^{5}+\frac{4056701}{406099094400}%
y^{6}+\frac{158419}{284269366080}y^{7} \\
&& \\
&&+\frac{11843}{304574320800}y^{10}+\frac{2461}{284269366080}y^{8}-\frac{%
158419}{1353663648}y^{4} \\
&& \\
&&+\frac{2461}{42301989}y^{3}-\frac{168263}{1705616196480}y^{9}-\frac{8999}{%
1705616196480}y^{11},
\end{eqnarray*}

\begin{eqnarray*}
d_{7}(y) &=&-\frac{2011}{6768318240}y^{5}-\frac{2011}{7614358020}y^{6}+\frac{%
158419}{1705616196480}y^{7} \\
&& \\
&&+\frac{11843}{1827445924800}y^{10}+\frac{2461}{1705616196480}y^{8}-\frac{%
158419}{8121981888}y^{4} \\
&& \\
&&+\frac{2461}{253811934}y^{3}-\frac{168263}{10233697178880}y^{9}-\frac{8999%
}{10233697178880}y^{11},
\end{eqnarray*}

\begin{eqnarray*}
d_{8}(y) &=&\frac{158419}{284269366080}y^{5}+\frac{158419}{1705616196480}%
y^{6}-\frac{19861}{596965668768}y^{7} \\
&& \\
&&-\frac{157481}{71635880252160}y^{10}+\frac{67}{682246478592}y^{8}+\frac{%
99305}{14213468304}y^{4} \\
&& \\
&&+\frac{335}{507623868}y^{3}+\frac{39253}{7163588025216}y^{9}+\frac{2687}{%
8954485031520}y^{11},
\end{eqnarray*}

\begin{eqnarray*}
d_{9}(y) &=&\frac{2461}{284269366080}y^{5}+\frac{2461}{1705616196480}y^{6}+%
\frac{67}{682246478592}y^{7} \\
&& \\
&&-\frac{8321}{71635880252160}y^{10}-\frac{1465}{2387862675072}y^{8}-\frac{%
335}{16243963776}y^{4} \\
&& \\
&&-\frac{7325}{1776683538}y^{3}+\frac{11251}{28654352100864}y^{9}+\frac{2003%
}{143271760504320}y^{11},
\end{eqnarray*}

\begin{eqnarray*}
d_{10}(y) &=&-\frac{168263}{1705616196480}y^{5}-\frac{168263}{10233697178880}%
y^{6}+\frac{39253}{7163588025216}y^{7} \\
&& \\
&&+\frac{38153}{85963056302592}y^{10}+\frac{11251}{28654352100864}y^{8}-%
\frac{196265}{170561619648}y^{4} \\
&& \\
&&+\frac{56255}{21320202456}y^{3}-\frac{6313}{5372691018912}y^{9}-\frac{12751%
}{214907640756480}y^{11},
\end{eqnarray*}

\begin{eqnarray*}
d_{11}(y) &=&\frac{11843}{304574320800}y^{5}+\frac{11843}{1827445924800}%
y^{6}-\frac{157481}{71635880252160}y^{7} \\
&& \\
&&-\frac{45611}{268634550945600}y^{10}-\frac{8321}{71635880252160}y^{8}+%
\frac{157481}{341123239296}y^{4} \\
&& \\
&&-\frac{8321}{10660101228}y^{3}+\frac{38153}{85963056302592}y^{9}+\frac{%
49001}{2149076407564800}y^{11},
\end{eqnarray*}

\begin{eqnarray*}
d_{12}(y) &=&-\frac{2687}{42640404912}y^{4}+\frac{2003}{21320202456}y^{3}-%
\frac{8999}{1705616196480}y^{5} \\
&& \\
&&-\frac{8999}{10233697178880}y^{6}+\frac{2687}{8954485031520}y^{7}+\frac{%
2003}{143271760504320}y^{8} \\
&& \\
&&-\frac{12751}{214907640756480}y^{9}-\frac{725}{236398404832128}y^{11}+%
\frac{49001}{2149076407564800}y^{10},
\end{eqnarray*}

\bigskip

\textbf{Proof:}%
\begin{eqnarray}
\left\langle u(x),{\normalsize R}_{y}{\normalsize (x)}\right\rangle
_{_{W_{2}^{6}}} &=&\sum_{i=0}^{5}u^{(i)}(0){\normalsize R}_{y}^{(i)}%
{\normalsize (0)}+\int_{0}^{1}u^{(6)}(x){\normalsize R}_{y}^{(6)}%
{\normalsize (x)}dx,\text{ }  \TCItag{2.3} \\
&&\left( u(x),\text{ }{\normalsize R}_{y}{\normalsize (x)}\in
W_{2}^{6}[0,1]\right) ,  \notag
\end{eqnarray}%
Through several integrations by parts for (2.3) we have

\begin{eqnarray}
\left\langle u(x),{\normalsize R}_{y}{\normalsize (x)}\right\rangle
_{_{W_{2}^{6}}} &=&\sum_{i=0}^{5}{\normalsize u}^{(i)}{\normalsize (0)}\left[
{\normalsize R}_{y}^{(i)}{\normalsize (0)-(-1)}^{(5-i)}{\normalsize R}%
_{y}^{(11-i)}{\normalsize (0)}\right]  \TCItag{2.4} \\
&&{\normalsize +}\sum_{i=0}^{5}{\normalsize (-1)}^{(5-i)}{\normalsize u}%
^{(i)}{\normalsize (1)R}_{y}^{(11-i)}{\normalsize (1)+}\int_{0}^{1}u%
{\normalsize (x)R}_{y}^{(12)}{\normalsize (x)dx.}  \notag
\end{eqnarray}%
Note that property of the reproducing kernel

\begin{equation*}
\left\langle u(x),\ R_{y}(x)\right\rangle _{W_{2}^{6}}=u(y),
\end{equation*}%
If 
\begin{equation}
\left\{ 
\begin{array}{c}
R_{y}^{(5)}(0)-R_{y}^{(6)}(0)=0, \\ 
R_{y}^{(4)}(0)+R_{y}^{(7)}(0)=0, \\ 
R_{y}^{\prime \prime \prime }(0)-R_{y}^{(8)}(0)=0, \\ 
R_{y}^{(6)}(1)=0, \\ 
R_{y}^{(7)}(1)=0, \\ 
R_{y}^{(8)}(1)=0, \\ 
R_{y}^{(9)}(1)=0,%
\end{array}%
\right.  \tag{2.5}
\end{equation}%
then (2.4) implies that,

\begin{equation*}
{\normalsize R}_{y}^{(12)}{\normalsize (x)=\delta (x-y),}
\end{equation*}%
When $x\neq y,$

\begin{equation*}
{\normalsize R}_{y}^{(12)}{\normalsize (x)=0,}
\end{equation*}%
therefore

\begin{equation*}
{\normalsize R}_{y}{\normalsize (x)=}\left\{ 
\begin{array}{c}
\sum_{i=1}^{12}\ c_{i}(y)x^{i-1},\text{ \ }x\leq y, \\ 
\\ 
\ \sum_{i=1}^{12}d_{i}(y)x^{i-1},\text{ \ }x>y.%
\end{array}%
\right.
\end{equation*}%
Since

\begin{equation*}
{\normalsize R}_{y}^{(12)}{\normalsize (x)=\delta (x-y),}
\end{equation*}%
we have

\begin{equation}
{\normalsize \partial }^{k}{\normalsize R}_{y^{+}}{\normalsize (y)=\partial }%
^{k}{\normalsize R}_{y^{-}}{\normalsize (y),}\text{ \ }%
k=0,1,2,3,4,5,6,7,8,9,10,  \tag{2.6}
\end{equation}%
and

\begin{equation}
{\normalsize \partial }^{11}{\normalsize R}_{y^{{\LARGE +}}}{\normalsize %
(y)-\partial }^{11}{\normalsize R}_{y^{-}}{\normalsize (y)=1.}  \tag{2.7}
\end{equation}%
Since ${\normalsize R}_{y}{\normalsize (x)\in }W_{2}^{6}[0,1]$, it follows
that

\begin{equation}
R_{y}(0)=0,\text{ \ }R_{y}(1)=0,\text{ }R_{y}^{\prime }(0)=0,\text{ \ }%
R_{y}^{\prime }(1)=0\text{, }R_{y}^{\prime \prime }(0)=0.  \tag{2.8}
\end{equation}%
From (2.5)-(2.8), the unknown coefficients $c_{i}(y)$ ve $d_{i}(y)$ $%
(i=1,2,...,12)$ can be obtained. Thus for $x\leq y,$ ${\normalsize R}_{y}%
{\normalsize (x)}$ is given by,

\begin{eqnarray*}
R_{y}(x) &=&\frac{2461}{42301989}x^{3}y^{5}+\frac{2461}{253811934}x^{3}y^{6}+%
\frac{335}{507623868}x^{3}y^{7}-\frac{8321}{10660101228}x^{3}y^{10} \\
&& \\
&&-\frac{7325}{1776683538}x^{3}y^{8}-\frac{11725}{84603978}x^{3}y^{4}+\frac{%
12221}{169207956}x^{3}y^{3}+\frac{56255}{21320202456}x^{3}y^{9} \\
&& \\
&&+\frac{2003}{21320202456}x^{3}y^{11}-\frac{158419}{1353663648}x^{4}y^{5}-%
\frac{158419}{8121981888}x^{4}y^{6}+\frac{99305}{14213468304}x^{4}y^{7} \\
&& \\
&&+\frac{157481}{341123239296}x^{4}y^{10}-\frac{335}{16243963776}x^{4}y^{8}+%
\frac{728021}{2707327296}x^{4}y^{4}-\frac{11725}{84603978}x^{4}y^{3} \\
&& \\
&&-\frac{196265}{170561619648}x^{4}y^{9}-\frac{2687}{42640404912}x^{4}y^{11}+%
\frac{4056701}{67683182400}x^{5}y^{5}-\frac{2011}{126905967}x^{5}y^{6} \\
&& \\
&&+\frac{158419}{284269366080}x^{5}y^{7}+\frac{11843}{304574320800}%
x^{5}y^{10}+\frac{2461}{284269366080}x^{5}y^{8}-\frac{158419}{135366364}%
x^{5}y^{4} \\
&& \\
&&+\frac{2461}{42301989}x^{5}y^{3}-\frac{168263}{1705616196480}x^{5}y^{9}-%
\frac{8999}{1705616196480}x^{5}y^{11}+\frac{4056701}{40609909440}x^{6}y^{5}
\\
&& \\
&&-\frac{2011}{7614358020}x^{6}y^{6}+\frac{158419}{1705616196480}x^{6}y^{7}+%
\frac{11843}{1827445924800}x^{6}y^{10}+\frac{2461}{170561619}x^{6}y^{8} \\
&& \\
&&-\frac{158419}{8121981888}x^{6}y^{4}+\frac{2461}{253811934}x^{6}y^{3}-%
\frac{168263}{10233697178880}x^{6}y^{9}-\frac{8999}{1023369717888}x^{6}y^{11}
\\
&& \\
&&+\frac{158419}{284269366080}x^{7}y^{5}+\frac{158419}{1705616196480}%
x^{7}y^{6}-\frac{19861}{596965668768}x^{7}y^{7}-\frac{157481}{716358802}%
x^{7}y^{10} \\
&& \\
&&+\frac{67}{682246478592}x^{7}y^{8}-\frac{104003}{81219818880}x^{7}y^{4}+%
\frac{335}{507623868}x^{7}y^{3}+\frac{39253}{7163588025216}x^{7}y^{9} \\
&& \\
&&+\frac{2687}{8954485031520}x^{7}y^{11}+\frac{2461}{284269366080}x^{8}y^{5}+%
\frac{2461}{1705616196480}x^{8}y^{6}+\frac{67}{682246478}x^{8}y^{7}
\end{eqnarray*}

\begin{eqnarray*}
&&-\frac{8321}{71635880252160}x^{8}y^{10}-\frac{1465}{2387862675072}%
x^{8}y^{8}-\frac{335}{16243963776}x^{8}y^{4}+\frac{12221}{1137077464320}%
x^{8}y^{3} \\
&& \\
&&+\frac{11251}{28654352100864}x^{8}y^{9}+\frac{2003}{143271760504320}%
x^{8}y^{11}-\frac{168263}{1705616196480}x^{9}y^{5}-\frac{168263}{1023369717}%
x^{9}y^{6} \\
&& \\
&&+\frac{39253}{7163588025216}x^{9}y^{7}+\frac{38153}{85963056302592}%
x^{9}y^{10}+\frac{11251}{28654352100864}x^{9}y^{8}-\frac{196265}{17056161964}%
x^{9}y^{4} \\
&& \\
&&+\frac{56255}{21320202456}x^{9}y^{3}-\frac{6313}{5372691018912}x^{9}y^{9}-%
\frac{12751}{214907640756480}x^{9}y^{11}-\frac{1}{725760}x^{9}y^{2} \\
&& \\
&&+\frac{11843}{304574320800}x^{10}y^{5}+\frac{11843}{1827445924800}%
x^{10}y^{6}-\frac{157481}{71635880252160}x^{10}y^{7}-\frac{45611}{2686345509}%
x^{10}y^{10} \\
&& \\
&&-\frac{8321}{71635880252160}x^{10}y^{8}+\frac{157481}{341123239296}%
x^{10}y^{4}-\frac{8321}{10660101228}x^{10}y^{3}+\frac{38153}{859630563025}%
x^{10}y^{9} \\
&& \\
&&+\frac{49001}{2149076407564800}y^{11}+\frac{1}{3628800}y-\frac{8999}{%
1705616196480}x^{11}y^{5}-\frac{8999}{10233697178880}x^{11}y^{6} \\
&& \\
&&+\frac{2687}{8954485031520}x^{11}y^{7}+\frac{49001}{2149076407564800}%
x^{11}y^{10}+\frac{2003}{143271760504320}x^{11}y^{8} \\
&& \\
&&-\frac{2687}{42640404912}x^{11}y^{4}+\frac{2003}{21320202456}x^{11}y^{3}-%
\frac{12751}{214907640756480}x^{11}y^{9}-\frac{725}{236398404832128}.
\end{eqnarray*}

\bigskip

\textbf{3. Solution represantation in }$W_{2}^{6}[0,1]$

In this section, the solution of equation (1.1) is given in the reproducing
kernel space $W_{2}^{6}[0,1].$

On defining the linear operator $L:W_{2}^{6}[0,1]\rightarrow W_{2}^{1}[0,1]$%
\ as

\begin{equation*}
Lu=u^{(5)}(x)-f(x)u(x).
\end{equation*}%
Model problem (1.1) changes the following problem:

\begin{equation}
\left\{ 
\begin{array}{c}
Lu=K(x),\text{ }x\in \lbrack 0,1] \\ 
u(a)=0,\text{ \ }u^{\prime }(a)=0,\text{ \ }u^{\prime \prime }(a)=0,\text{ \ 
}u(b)=0,\text{ \ }u^{\prime }(b)=0.\text{\ }%
\end{array}%
\right. \text{\ }  \tag{3.1}
\end{equation}

\bigskip

\textbf{3.1. The Linear boundedness of operator }$L.$

\textbf{Lemma 3.1. }If $u(x)\in W_{2}^{6}[a,b],$ then $\left\Vert
u^{(k)}(x)\right\Vert _{L^{\infty }}\leq M_{k}\left\Vert u(x)\right\Vert
_{W_{2}^{6}},$ where $M_{k}$ $(k=0,1,...,5)~$are positive constants.

\textbf{Proof: }For any $x\in \lbrack a,b]$ it holds that

\begin{equation*}
\left\Vert R_{x}(y)\right\Vert _{W_{2}^{6}}=\sqrt{\left\langle
R_{x}(y),R_{x}(y)\right\rangle _{W_{2}^{6}}}=\sqrt{R_{x}(x),}
\end{equation*}%
from the continuity of $R_{x}(x),$ there exists a constant $M_{0}>0,$ such
that $\left\Vert R_{x}(y)\right\Vert _{W_{2}^{6}}\leq M_{0}.$ By (2.1) one
gets

\begin{equation}
\left\vert u(x)\right\vert =\left\vert \left\langle
u(y),R_{x}(y)\right\rangle _{_{W_{2}^{6}}}\right\vert \leq \left\Vert
R_{x}(y)\right\Vert _{W_{2}^{6}}\left\Vert u(y)\right\Vert _{W_{2}^{6}}\leq
M_{0}\left\Vert u(y)\right\Vert _{W_{2}^{6}}.  \tag{3.2}
\end{equation}%
For any $x,$ $y\in \lbrack a,b]$, there exists $M_{k}$ $(k=1,2,...,5),$ such
that

\begin{equation*}
\left\Vert R_{x}^{(k)}(y)\right\Vert _{W_{2}^{6}}\leq M_{k}\text{ \ }%
(k=1,2,...,5),
\end{equation*}%
we have%
\begin{eqnarray}
\left\vert u^{(k)}(x)\right\vert &=&\left\vert \left\langle
u(y),R_{x}^{(k)}(y)\right\rangle _{_{W_{2}^{6}}}\right\vert \leq \left\Vert
R_{x}^{(k)}(y)\right\Vert _{W_{2}^{6}}\left\Vert u(y)\right\Vert
_{W_{2}^{6}}\leq M_{k}\left\Vert u(y)\right\Vert _{W_{2}^{6}}\text{ } 
\TCItag{3.3} \\
(k &=&1,2,...,5).  \notag
\end{eqnarray}%
Combining (3.2) and (3.3), it follows that

\begin{equation*}
\left\Vert u^{(k)}(x)\right\Vert _{L^{\infty }}\leq M_{k}\left\Vert
u(y)\right\Vert _{W_{2}^{6}}\text{ }(k=1,2,...,5).
\end{equation*}

\bigskip

\textbf{Theorem 3.1. }Suppose $f_{i}^{\prime }\in L^{2}[a,b]$ $%
(i=0,1,2,3,4), $ Then $L:W_{2}^{6}[a,b]\rightarrow W_{2}^{1}[a,b]$ is a
bounded linear operator.

\textbf{Proof: }(i)\textbf{\ }By the definition of the operator it is clear
that $L$ is a linear operator.

(ii) Due to the definiton of $W_{2}^{1}[a,b]$ , we have%
\begin{eqnarray*}
\left\Vert (Lu)(x)\right\Vert _{W_{2}^{1}}^{2} &=&\left\langle
(Lu)(x),(Lu)(x)\right\rangle _{W_{2}^{1}} \\
&=&[(Lu)(a)]^{2}+\int_{a}^{b}[(Lu)^{\prime }(x)]^{2}dx \\
&=&\left[ \sum_{i=0}^{5}f_{i}(a)u^{(i)}(a)\right] ^{2}+\int_{a}^{b}\left[
\left( \sum_{i=0}^{5}f_{i}(x)u^{(i)}(x)\right) ^{\prime }\right] ^{2}dx.
\end{eqnarray*}

\begin{eqnarray*}
\int_{a}^{b}[(Lu)^{\prime }(x)]^{2}dx &=&\int_{a}^{b}\left[
u^{(6)}(x)+\sum_{i=0}^{4}(f_{i}^{\prime }(x)u^{(i)}(x)+f_{i}(x)u^{(i+1)}(x))%
\right] ^{2}dx \\
&=&\int_{a}^{b}\left[ u^{(6)}(x)\right] ^{2}dx+2\int_{a}^{b}\left[
u^{(6)}(x)\sum_{i=0}^{4}(f_{i}^{\prime }(x)u^{(i)}(x)+f_{i}(x)u^{(i+1)}(x))%
\right] dx \\
&&+\int_{a}^{b}\left[ \sum_{i=0}^{4}(f_{i}^{\prime
}(x)u^{(i)}(x)+f_{i}(x)u^{(i+1)}(x))\right] ^{2}dx,
\end{eqnarray*}%
where

\begin{equation*}
\int_{a}^{b}\left[ u^{(6)}(x)\right] ^{2}dx\leq \left\Vert u(x)\right\Vert
_{W_{2}^{6}}^{2},
\end{equation*}%
and

\begin{eqnarray*}
&&\int_{a}^{b}\left[ u^{(6)}(x)\sum_{i=0}^{4}(f_{i}^{\prime
}(x)u^{(i)}(x)+f_{i}(x)u^{(i+1)}(x))\right] dx \\
&\leq &\left\{ \int_{a}^{b}\left[ u^{(6)}(x)\right] ^{2}dx\right\}
^{2}\left\{ \int_{a}^{b}\left[ \sum_{i=0}^{4}(f_{i}^{\prime
}(x)u^{(i)}(x)+f_{i}(x)u^{(i+1)}(x))\right] ^{2}dx\right\} ^{2}.
\end{eqnarray*}%
By lemma $\left( 3.1\right) $ and $f_{i}^{\prime }(x)$ $\in L^{2}[a,b],$ we
can obtain a constant $N>0,$ satisfying

\begin{equation*}
\int_{a}^{b}\left[ \sum_{i=0}^{4}(f_{i}^{\prime
}(x)u^{(i)}(x)+f_{i}(x)u^{(i+1)}(x))\right] ^{2}dx\leq N(b-a)\left\Vert
u(x)\right\Vert _{W_{2}^{6}}^{2}.
\end{equation*}%
Furthermore one gets

\begin{equation*}
\int_{a}^{b}[(Lu)^{\prime }(x)]^{2}dx\leq \left\Vert u(x)\right\Vert
_{W_{2}^{6}}^{2}+2\sqrt{N(b-a)}\left\Vert u(x)\right\Vert
_{W_{2}^{6}}^{2}+N(b-a)\left\Vert u(x)\right\Vert _{W_{2}^{6}}^{2},
\end{equation*}%
let $G=\left( 1+2\sqrt{N(b-a)}+N(b-a)\right) >0,$ then

\begin{equation*}
\int_{a}^{b}[(Lu)^{\prime }(x)]^{2}dx\leq G\left\Vert u(x)\right\Vert
_{W_{2}^{6}}^{2}.
\end{equation*}%
Therefore $L$ is a boundend operator. So we obtain the result as required. \
\ \ \ \ $\square $

\bigskip

\textbf{3.2. The normal orthogonal function system of }$W_{2}^{6}[a,b]$%
\bigskip

We choose $\left\{ x_{i}\right\} _{i=1}^{\infty }$ as any dense set in $%
[a,b] $ and let $\Psi _{x}(y)=L^{\ast }T_{x}(y),$ where $L^{\ast \text{ }}$
is conjugate operator of $L$ and $T_{x}(y)$ is given by (2.2). Furthermore,
for simplicity let $\Psi _{i}(x)=\Psi _{x_{i}}(x),$ namely,

\begin{equation*}
\Psi _{i}(x)\overset{def}{=}\Psi _{x_{i}}(x)=\text{ }L^{\ast \text{ }%
}T_{x_{i}}(x).
\end{equation*}%
Now several lemmas are given.

\bigskip

\textbf{Lemma 3.2. }$\left\{ \Psi _{i}(x)\right\} _{i=1}^{\infty }$ is
complete system of $W_{2}^{6}[a,b].$\bigskip

\textbf{Proof: }For $u(x)\in W_{2}^{6}[a,b]$, let $\left\langle u(x),\Psi
_{i}(x)\right\rangle =0$ $(i=1,2,...),$ that is

\begin{equation*}
\left\langle u(x),\text{ }L^{\ast \text{ }}T_{x_{i}}(x)\right\rangle
=(Lu)(x_{i})=0.
\end{equation*}%
Note that $\left\{ x_{i}\right\} _{i=1}^{\infty }$ is the dense set in $%
[a,b],$ therefore $(Lu)(x)=0.$ It follows that $u(x)=0$ from the existence
of $L^{-1}.$ \ \ \ \ $\square $

\bigskip 

\textbf{Lemma 3.3. }The following formula holds

\begin{equation*}
\Psi _{i}(x)=\left( L\eta R_{x}(\eta )\right) \left( x_{i}\right) ,
\end{equation*}%
where the subscript $\eta $ of operator $L\eta $ indicates that the operator 
$L$ applies to function of $\eta .$

\bigskip

\textbf{Proof:}%
\begin{eqnarray*}
\Psi _{i}(x) &=&\left\langle \Psi _{i}(\xi ),R_{x}(\xi )\right\rangle
_{W_{2}^{6}[a,b]} \\
&=&\left\langle L^{\ast \text{ }}T_{x_{i}}\left( \xi \right) ,R_{x}(\xi
)\right\rangle _{W_{2}^{6}[a,b]} \\
&=&\left\langle \left( T_{x_{i}}\right) \left( \xi \right) ,\left( L_{\eta
}R_{x}(\eta )\right) \left( \xi \right) \right\rangle _{W_{2}^{1}[a,b]} \\
&=&\left( L_{\eta }R_{x}(\eta )\right) \left( x_{i}\right) .
\end{eqnarray*}%
This completes the proof.\bigskip\ \ \ \ \ \ \ \ $\square $

\textbf{Remark 3.1.} The orthonormal system $\left\{ \overline{\Psi }%
_{i}(x)\right\} _{i=1}^{\infty }$ of $W_{2}^{6}[a,b]$ can be derived from
Gram-Schmidt orthogonaliztion process of $\left\{ \Psi _{i}(x)\right\}
_{i=1}^{\infty },$

\begin{equation}
\overline{\Psi }_{i}(x)=\sum_{k=1}^{i}\beta _{ik}\Psi _{k}(x),\text{ \ }%
(\beta _{ii}>0,\text{ \ }i=1,2,...)  \tag{3.7}
\end{equation}%
where $\beta _{ik}$ are orthogonal cofficients.\bigskip

In the following, we will give the represantation of the exact solution of
Eq.(1.1) in the reproducing kernel space $W_{2}^{6}[a,b].$\bigskip

\textbf{3.3. The structure of the solution and the main results}\bigskip

\textbf{Theorem 3.2.} If $u(x)$ is the exact solution of Eq.(1.1), then

\begin{equation*}
u(x)=\sum_{i=1}^{\infty }\sum_{k=1}^{i}\beta _{ik}g(x_{k})\overline{\Psi }%
_{i}(x),
\end{equation*}%
where $\left\{ x_{i}\right\} _{i=1}^{\infty }$ is a dense set in $[a,b].$%
\bigskip

\textbf{Proof: }From the (3.7) and uniqeness of solution of (1.1) (see
[32]), we have

\begin{eqnarray*}
u(x) &=&\sum_{i=1}^{\infty }\left\langle u(x),\overline{\Psi }%
_{i}(x)\right\rangle _{W_{2}^{6}}\overline{\Psi }_{i}(x) \\
&=&\sum_{i=1}^{\infty }\sum_{k=1}^{i}\beta _{ik}\left\langle u(x),L^{\ast
}T_{x_{k}}(x)\right\rangle _{W_{2}^{6}}\overline{\Psi }_{i}(x) \\
&=&\sum_{i=1}^{\infty }\sum_{k=1}^{i}\beta _{ik}\left\langle
Lu(x),T_{x_{k}}(x)\right\rangle _{W_{2}^{1}}\overline{\Psi }_{i}(x) \\
&=&\sum_{i=1}^{\infty }\sum_{k=1}^{i}\beta _{ik}\left\langle
g(x),T_{x_{k}}(x)\right\rangle _{W_{2}^{1}}\overline{\Psi }_{i}(x) \\
&=&\sum_{i=1}^{\infty }\sum_{k=1}^{i}\beta _{ik}g(x_{k})\overline{\Psi }%
_{i}(x).\text{ \ \ \ \ \ \ }\square
\end{eqnarray*}%
Now the approximate solution $u_{n}(x)$ can be obtained by truncating the $%
n- $ term of the exact solution $u(x),$

\begin{equation*}
u_{n}(x)=\sum_{i=1}^{n}\sum_{k=1}^{i}\beta _{ik}g(x_{k})\overline{\Psi }%
_{i}(x).
\end{equation*}

\textbf{Theorem 3.3. }Assume $u(x)$ is the solution of Eq.(1.1) and $r_{n}(x)
$ is the error between the approximate solution $u_{n}(x)$ and the exact
solution $u(x).$ Then the error sequence $r_{n}(x)$ is monotone decreasing
in the sense of $\left\Vert .\right\Vert _{W_{2}^{6}}$ and $\left\Vert
r_{n}(x)\right\Vert _{W_{2}^{6}}\rightarrow 0$ [23].

\bigskip 

\textbf{4. Numerical Results}

In this section, four numerical examples are provided to show the accuracy
of the present method. All computations are performed by Maple 13. The RKHSM
does not require discretization of the variables, i.e., time and space, it
is not effected by computation round off errors and one is not faced with
necessity of large computer memory and time. The accuracy of the RKHSM for
the fifth-order boundary value problems is controllable and absolute errors
are small with present choice of $x$ (see Table 1-4). The numerical results
we obtained justify the advantage of this methodology.\bigskip 

\textbf{Example 4.1.} ([8,11,12]). We first consider the linear boundary
value problem%
\begin{equation}
\left\{ 
\begin{array}{c}
y^{(5)}(x)=y-15e^{x}-10xe^{x},\text{ }0<x<1 \\ 
y(0)=0,\text{ }y^{\prime }(0)=1,\text{ }y^{\prime \prime }(0)=0,\text{ }%
y(1)=0,\text{ }y^{\prime }(1)=-e%
\end{array}%
\right.  \tag{4.1}
\end{equation}%
The exact solution of (4.1) is

\begin{equation*}
y\left( x\right) =x(1-x)e^{x}.
\end{equation*}

If we homogenize the boundary conditions of (4.1), then the following (4.2)
is obtained%
\begin{equation}
\left\{ 
\begin{array}{c}
u^{(5)}(x)-u(x)=1-5e^{x}\left[ 2-3x+3x^{2}(2-\dfrac{5}{e})+4x^{3}(\dfrac{-3}{%
2}+\dfrac{4}{e})\right] \\ 
\\ 
-10e^{x}\left[ -3+6x(2-\dfrac{5}{e})+12x^{2}(\dfrac{-3}{2}+\dfrac{4}{e})%
\right] \\ 
\\ 
-10e^{x}\left[ 6(2-\dfrac{5}{e})+24x(\dfrac{-3}{2}+\dfrac{4}{e})\right]
-5e^{x}\left[ 24(\dfrac{-3}{2}+\dfrac{4}{e})\right] \\ 
\\ 
-15e^{x}-10xe^{x},\text{ \ \ \ \ \ \ \ \ \ \ \ \ \ \ \ }0<x<1 \\ 
\\ 
u(0)=0,\text{ }u^{\prime }(0)=0,\text{ }u^{\prime \prime }(0)=0,\text{ }%
u(1)=0,\text{ }u^{\prime }(1)=0%
\end{array}%
\right.  \tag{4.2}
\end{equation}%
\bigskip

\textbf{Example 4.2.} ([8,12]). We now consider the nonlinear BVP%
\begin{equation}
\left\{ 
\begin{array}{c}
y^{(5)}(x)=e^{-x}y^{2}(x),\text{ }0<x<1 \\ 
y(0)=1,\text{ }y^{\prime }(0)=1,\text{ }y^{\prime \prime }(0)=1,\text{ }%
y(1)=e,\text{ }y^{\prime }(1)=e%
\end{array}%
\right.  \tag{4.3}
\end{equation}%
The exact solution of (4.3) is

\begin{equation*}
y(x)=e^{x}
\end{equation*}

If we homogenize the boundary conditions of (4.3), then the following (4.4)
is obtained%
\begin{equation}
\left\{ 
\begin{array}{c}
u^{(5)}(x)-2e^{-x}\left[ 1+e^{x}(x-\dfrac{x^{2}}{2}+x^{3}(2-\dfrac{5}{e}%
)+x^{4}(\dfrac{-3}{2}+\dfrac{4}{e}))\right] u(x) \\ 
\\ 
=e^{-x}u^{2}(x)+e^{-x}\left[ 1+e^{x}(x-\dfrac{x^{2}}{2}+x^{3}(2-\dfrac{5}{e}%
)+x^{4}(\dfrac{-3}{2}+\dfrac{4}{e}))\right] ^{2} \\ 
\\ 
-e^{x}\left[ x-\dfrac{x^{2}}{2}+x^{3}(2-\dfrac{5}{e})+x^{4}(\dfrac{-3}{2}+%
\dfrac{4}{e})\right] \\ 
\\ 
-5e^{x}\left[ 1-x+3x^{2}(2-\dfrac{5}{e})+4x^{3}(\dfrac{-3}{2}+\dfrac{4}{e})%
\right] \\ 
\\ 
-10e^{x}\left[ -1+6x(2-\dfrac{5}{e})+12x^{2}(\dfrac{-3}{2}+\dfrac{4}{e})%
\right] \\ 
\\ 
-10e^{x}\left[ 6(2-\dfrac{5}{e})+24x(\dfrac{-3}{2}+\dfrac{4}{e})\right] \\ 
\\ 
-5e^{x}\left[ 24(\dfrac{-3}{2}+\dfrac{4}{e})\right] ,\text{ \ \ \ \ \ \ }%
0<x<1\text{\ \ \ \ \ \ \ \ } \\ 
\\ 
u(0)=0,\text{ }u^{\prime }(0)=0,\text{ }u^{\prime \prime }(0)=0,\text{ }%
u(1)=0,\text{ }u^{\prime }(1)=0%
\end{array}%
\right.  \tag{4.4}
\end{equation}

\bigskip \textbf{Example 4.3.([15]). }Consider the nonlinear BVP%
\begin{equation}
\left\{ 
\begin{array}{c}
y^{(5)}(x)=-24e^{-y(x)}+\frac{{\LARGE 48}}{{\LARGE (1+x)}^{5}},\text{ \ \ \ }%
0<x<1 \\ 
y(0)=0,\text{ }y^{\prime }(0)=1,\text{ }y^{\prime \prime }(0)=-1,\text{ }%
y(1)=\ln 2,\text{ }y^{\prime }(1)=0.5%
\end{array}%
\right. .  \tag{4.5}
\end{equation}%
The exact solution of (4.5) is

\begin{equation*}
y\left( x\right) =\ln (x+1).
\end{equation*}

If we homogenize the boundary conditions of (4.5), then the following (4.6)
is obtained

\begin{equation}
\left\{ 
\begin{array}{c}
u^{(5)}(x)=-24e^{-\left( {\LARGE u(x)}+{\LARGE x}-\dfrac{x^{2}}{2}+{\LARGE x}%
^{3}({\LARGE 4\ln 2-}\frac{{\LARGE 5}}{{\LARGE 2}}{\LARGE )}+{\LARGE x}^{4}%
{\LARGE (2-3}\ln {\LARGE 2)}\right) }+\frac{{\LARGE 48}}{{\LARGE (1+x)}^{5}}
\\ 
u(0)=0,\text{ }u^{\prime }(0)=0,\text{ }u^{\prime \prime }(0)=0,\text{ }%
u(1)=0,\text{ }u^{\prime }(1)=0%
\end{array}%
\right.  \tag{4.6}
\end{equation}

\bigskip

\textbf{Example 4.4. ([15]).} This is the nonlinear BVP%
\begin{equation}
\left\{ 
\begin{array}{c}
y^{(5)}(x)+y^{(4)}(x)+e^{-2x}y^{2}(x)=2e^{x}+1\text{ \ \ \ }0<x<1 \\ 
y(0)=0,\text{ }y^{\prime }(0)=1,\text{ }y^{\prime \prime }(0)=1,\text{ }%
y(1)=e,\text{ }y^{\prime }(1)=e%
\end{array}%
\right. .  \tag{4.7}
\end{equation}%
The exact solution of (4.7) is

\begin{equation*}
y\left( x\right) =e^{x}.
\end{equation*}

If we homogenize the boundary conditions of (4.7), then the following (4.8)
is obtained

\begin{equation}
\left\{ 
\begin{array}{c}
u^{(5)}(x)+u^{(4)}(x)=-e^{-2x}(u(x)+1+x+\frac{x^{2}}{2}+x^{3}(3e-8)+x^{4}(%
\frac{11}{2}-2e))^{2}+2e^{x}+48e-131 \\ 
u(0)=0,\text{ }u^{\prime }(0)=0,\text{ }u^{\prime \prime }(0)=0,\text{ }%
u(1)=0,\text{ }u^{\prime }(1)=0%
\end{array}%
\right.  \tag{4.8}
\end{equation}

\bigskip

\textbf{Remark 4.1.} Lamnii et al. [16] solved the problem (5.1) by using
sextic spline collocation method. He obtained the accurate approximate
solutions of this problem for the small $h$ values. Zhang [44] investigated
approximate solution of the problem (5.1) by using variational iteration
method. In addition, the same problem is solved by Noor and Mohyud-Din [12]
previously and they got better results by using the variational iteration
method.

Lv and Cui [32] studied only the linear fifth-order two-point boundary value
problems by using reproducing kernel Hilbert space method. We use the RKHSM
for the same linear problem with different boundary conditions and also use
different reproducing kernel function for computations. We also use the
RKHSM for the nonlinear problems.

Using our method we chose $36$ points on $[0,1].$ In Tables 1-4, we computed
the absolute errors $\left\vert u\left( x,t\right) -u_{n}\left( x,t\right)
\right\vert $ at the points $\left\{ \left( x_{i}\right) :\text{ }x_{i}=i,%
\text{ \ }i=0.0,0.1,...,1.0\right\} .$

\bigskip

\begin{center}
\begin{equation*}
\begin{tabular}{|l|l|l|l|l|}
\hline
$x$ & Exact Sol. & RKHSM & HPM [13] & B-Spline [11] \\ \hline
$0.0$ & $0.0000$ & $0.000$ & $0.0000$ & $0.0000$ \\ \hline
$0.1$ & $0.099465382$ & $5.89\times 10^{-7}$ & $3\times 10^{-11}$ & $%
8.0\times 10^{-3}$ \\ \hline
$0.2$ & $0.195424441$ & $1.73\times 10^{-8}$ & $2\times 10^{-10}$ & $%
1.2\times 10^{-3}$ \\ \hline
$0.3$ & $0.283470349$ & $6.02\times 10^{-7}$ & $4\times 10^{-10}$ & $%
5.0\times 10^{-3}$ \\ \hline
$0.4$ & $0.358037927$ & $7.42\times 10^{-7}$ & $8\times 10^{-10}$ & $%
3.0\times 10^{-3}$ \\ \hline
$0.5$ & $0.412180317$ & $3.32\times 10^{-7}$ & $1.2\times 10^{-9}$ & $%
8.0\times 10^{-3}$ \\ \hline
$0.6$ & $0.437308512$ & $3.10\times 10^{-7}$ & $2\times 10^{-9}$ & $%
6.0\times 10^{-3}$ \\ \hline
$0.7$ & $0.422888068$ & $3.08\times 10^{-7}$ & $2.2\times 10^{-9}$ & $0.000$
\\ \hline
$0.8$ & $0.356086548$ & $4.58\times 10^{-7}$ & $1.9\times 10^{-9}$ & $%
9.0\times 10^{-3}$ \\ \hline
$0.9$ & $0.221364280$ & $4.30\times 10^{-7}$ & $1.4\times 10^{-9}$ & $%
9.0\times 10^{-3}$ \\ \hline
$1.0$ & $0.0000$ & $2.36\times 10^{-13}$ & $0.0000$ & $0.0000$ \\ \hline
\end{tabular}%
\begin{tabular}{|l|l|}
\hline
ADM [8] & Sinc [15] \\ \hline
$0.0000$ & $0.0000$ \\ \hline
$3\times 10^{-11}$ & $0.0000$ \\ \hline
$2\times 10^{-10}$ & $0.1\times 10^{-5}$ \\ \hline
$4\times 10^{-10}$ & $0.3\times 10^{-5}$ \\ \hline
$8\times 10^{-10}$ & $0.3\times 10^{-5}$ \\ \hline
$1.2\times 10^{-9}$ & $0.0000$ \\ \hline
$2\times 10^{-9}$ & $0.5\times 10^{-5}$ \\ \hline
$2.2\times 10^{-9}$ & $0.9\times 10^{-5}$ \\ \hline
$1.9\times 10^{-9}$ & $0.2\times 10^{-5}$ \\ \hline
$1.4\times 10^{-9}$ & $0.1\times 10^{-5}$ \\ \hline
$0.0000$ & $\ 0.0000$ \\ \hline
\end{tabular}%
\end{equation*}%
\\[0pt]

\textbf{Table 1.} The absolute error of Example 5.1 for boundary conditions
at $0.0\leq x\leq 1.0.$

\bigskip

\begin{equation*}
\begin{tabular}{|l|l|l|l|l|}
\hline
$x$ & Exact Solution & RKHSM & HPM [13] & B-Spline[11] \\ \hline
$0.0$ & $0.$ & $0.0000$ & $0.0000$ & $0.0000$ \\ \hline
$0.1$ & $1.105170918$ & $5.19\times 10^{-7}$ & $1\times 10^{-9}$ & $7.0$ $%
\times 10^{-4}$ \\ \hline
$0.2$ & $1.221402758$ & $0.60\times 10^{-7}$ & $2\times 10^{-9}$ & $7.2$ $%
\times 10^{-4}$ \\ \hline
$0.3$ & $1.349858808$ & $3.19\times 10^{-7}$ & $1\times 10^{-9}$ & $4.1$ $%
\times 10^{-4}$ \\ \hline
$0.4$ & $1.491824698$ & $2.50\times 10^{-7}$ & $2\times 10^{-8}$ & $4.6$ $%
\times 10^{-4}$ \\ \hline
$0.5$ & $1.648721271$ & $3.03\times 10^{-7}$ & $3.1\times 10^{-8}$ & $4.7$ $%
\times 10^{-4}$ \\ \hline
$0.6$ & $1.822118800$ & $9.60\times 10^{-7}$ & $3.7\times 10^{-8}$ & $4.8$ $%
\times 10^{-4}$ \\ \hline
$0.7$ & $2.013752707$ & $4.20\times 10^{-7}$ & $4.1\times 10^{-8}$ & $3.9$ $%
\times 10^{-4}$ \\ \hline
$0.8$ & $2.225540928$ & $4.09\times 10^{-7}$ & $3.1\times 10^{-8}$ & $3.1$ $%
\times 10^{-4}$ \\ \hline
$0.9$ & $2.459603111$ & $5.46\times 10^{-7}$ & $1.4\times 10^{-8}$ & $1.6$ $%
\times 10^{-4}$ \\ \hline
$1.0$ & $2.718281828$ & $5.34\times 10^{-7}$ & $0.0000$ & $0.0000$ \\ \hline
\end{tabular}%
\begin{tabular}{|l|l|}
\hline
ADM [8] & VIM [12] \\ \hline
\multicolumn{1}{|l|}{$0.0000$} & $0.0000$ \\ \hline
\multicolumn{1}{|l|}{$1\times 10^{-9}$} & $1\times 10^{-9}$ \\ \hline
\multicolumn{1}{|l|}{$2\times 10^{-9}$} & $2\times 10^{-9}$ \\ \hline
\multicolumn{1}{|l|}{$1\times 10^{-9}$} & $1\times 10^{-9}$ \\ \hline
\multicolumn{1}{|l|}{$2\times 10^{-8}$} & $2\times 10^{-8}$ \\ \hline
\multicolumn{1}{|l|}{$3.1\times 10^{-8}$} & $3.1\times 10^{-8}$ \\ \hline
\multicolumn{1}{|l|}{$3.7\times 10^{-8}$} & $3.7\times 10^{-8}$ \\ \hline
\multicolumn{1}{|l|}{$4.1\times 10^{-8}$} & $4.1\times 10^{-8}$ \\ \hline
\multicolumn{1}{|l|}{$3.1\times 10^{-8}$} & $3.1\times 10^{-8}$ \\ \hline
\multicolumn{1}{|l|}{$1.4\times 10^{-8}$} & $1.4\times 10^{-8}$ \\ \hline
\multicolumn{1}{|l|}{$0.0000$} & $0.0000$ \\ \hline
\end{tabular}%
\end{equation*}

\textbf{Table 2.} The absolute error of Example 5.2 for boundary conditions
at $0.0\leq x\leq 1.0.$

\bigskip
\end{center}

\begin{equation*}
\begin{tabular}{|l|l|l|l|}
\hline
x & Exact Solution & RKHSM & 
\begin{tabular}{l}
AE, $1.0E-8$ \\ 
RKHSM%
\end{tabular}
\\ \hline
$0.0$ & $0.0$ & $0.0$ & $0.0$ \\ \hline
$0.0806$ & $0.07751644243$ & $0.07751634304$ & $0.003$ \\ \hline
$0.1648$ & $0.1525493985$ & $0.1525493860$ & $1.25$ \\ \hline
$0.2285$ & $0.2057939130$ & $0.2057939053$ & $0.77$ \\ \hline
$0.3999$ & $0.3364008055$ & $0.3364007134$ & $9.21$ \\ \hline
$0.5$ & $0.4054651081$ & $0.4054650667$ & $4.14$ \\ \hline
$0.6923$ & $0.5260885504$ & $0.5260885142$ & $3.62$ \\ \hline
$0.7714$ & $0.5717701944$ & $0.5717701752$ & $1.92$ \\ \hline
$0.8836$ & $0.6331848394$ & $0.6331848843$ & $4.49$ \\ \hline
$0.9447$ & $0.6651077235$ & $0.6651077287$ & $0.52$ \\ \hline
$1.0$ & $0.6931471806$ & $0.6931471783$ & $0.23$ \\ \hline
\end{tabular}%
\begin{tabular}{|l|}
\hline
\begin{tabular}{l}
AE,$1.0E-4$ [15] \\ 
Sinc-Galerkin%
\end{tabular}
\\ \hline
$0.0$ \\ \hline
$0.0$ \\ \hline
$0.2$ \\ \hline
$0.2$ \\ \hline
$0.4$ \\ \hline
$0.1$ \\ \hline
$0.2$ \\ \hline
$0.3$ \\ \hline
$0.2$ \\ \hline
$0.5$ \\ \hline
$0.0$ \\ \hline
\end{tabular}%
\end{equation*}

\begin{center}
\textbf{Table 3.} The absolute error (AE) of Example 5.3 for boundary
conditions at $0.0\leq x\leq 1.0.$

\bigskip
\end{center}

\begin{equation*}
\begin{tabular}{|l|l|l|l|}
\hline
x & Exact Solution & RKHSM & 
\begin{tabular}{l}
AE \\ 
RKHSM%
\end{tabular}
\\ \hline
$0.0$ & $1.0$ & $1.0$ & $0.0$ \\ \hline
$0.0100$ & $1.105170918$ & $1.105170918$ & $0.0$ \\ \hline
$0.1184$ & $1.125694299$ & $1.125694299$ & $0.0$ \\ \hline
$0.1517$ & $1.163811041$ & $1.163811041$ & $0.0$ \\ \hline
$0.2410$ & $1.272521035$ & $1.272521035$ & $0.0$ \\ \hline
$0.3604$ & $1.433902861$ & $1.433902861$ & $0.0$ \\ \hline
$0.4287$ & $1.535260387$ & $1.535260387$ & $0.0$ \\ \hline
$0.5000$ & $1.648721271$ & $1.648721271$ & $0.0$ \\ \hline
$0.6395$ & $1.895532876$ & $1.895532876$ & $0.0$ \\ \hline
$0.8482$ & $2.335439276$ & $2.335439276$ & $0.0$ \\ \hline
$0.9996$ & $2.717194733$ & $2.717194734$ & $1\times 10^{-9}$ \\ \hline
$1.0$ & $2.718281828$ & $2.718281828$ & $0.0$ \\ \hline
\end{tabular}%
\begin{tabular}{|l|}
\hline
\begin{tabular}{l}
AE,$1.0E-3$ [15] \\ 
Sinc-Galerkin%
\end{tabular}
\\ \hline
$0.0$ \\ \hline
$0.0$ \\ \hline
$0.0$ \\ \hline
$0.1$ \\ \hline
$0.0$ \\ \hline
$0.1$ \\ \hline
$0.0$ \\ \hline
$0.2$ \\ \hline
$0.1$ \\ \hline
$0.2$ \\ \hline
$0.2$ \\ \hline
$0.0$ \\ \hline
\end{tabular}%
\end{equation*}

\begin{center}
\textbf{Table 4.} The absolute error of Example 4 for boundary conditions at 
$0.0\leq x\leq 1.0.$

\bigskip
\end{center}

\textbf{\ Remark 4.2.} The RKHSM tested on four problems, one linear and
three nonlinears. A comparison with decomposition method by Wazwaz [8],
sixth B-spline method by Caglar et al. [11], variational iteration and
homotopy perturbation methods by Noor and Mohyid-Din [12,13] and
sinc-galerkin method by Gamel [15] are made and it was seen that the present
method yields good results (see Tables 1-4).\bigskip

\textbf{6. Conclusion}

In this paper, we introduce an algorithm for solving the fifth-order problem
with boundary conditions. For illustration purposes, chose four examples
which were selected to show the computational accuracy. It may be concluded
that, the RKHSM is very powerful and efficient in finding exact solution for
a wide class of boundary value problems. The method gives more realistic
series solutions that converge very rapidly in physical problems. The
approximate solution obtained by the present method is uniformly convergent.

Clearly, the series solution methodology can be applied to much more
complicated nonlinear differential equations and boundary value problems
[20-43] as well. However, if the problem becomes nonlinear, then the RKHSM
does not require discretization or perturbation and it does not make closure
approximation. Results of numerical examples show that the present method is
an accurate and reliable analytical method for the fifth order problem with
boundary conditions.

\bigskip

\textbf{References}

[1] M. R. Scott, H. A. Watts, SUPORT-A computer code for two-point
boundary-value problems via orthonormalization, SAND75-0198, Sandia
Laboratories, Albuquerque, NM, 1975.

[2] M. R. Scott, H. A. Watts, Computational solution of linear two-point
boundary value problems via orthonormalization, SIAM J. Numer. Anal. 14 (1)
(1977) 40--70.

[3] M. R. Scott, W. H. Vandevender, A comparison of several invariant
imbedding algorithms for the solution of two-point boundary-value problems,
Appl. Math. Comput. 1 (1975) 187--218.

[4] A. R. Davies, A. Karageorghis, T.N. Phillips, Spectral Galerkin method
for primary two-point boundary-value problems in modelling viscoelastic
flows, Int. J. Num. Methods Eng. 26 (1988) 647-662.

[5] A. Karageorghis, T. N. Phillips, A.R. Davies, Spectral collocation
methods for the primary two-point boundary-value problems in modelling
viscoelastic flows, Int. J. Num. Methods Eng. 26 (1988) 805-813.

[6] R. P. Agarwal, Boundary Value Problems for High Order Differential
Equations, World Scientific, Singapore, 1986.

[7] M. S. Khan, Finite-difference solutions of fifth-order boundary-value
problems, Ph.D. Thesis, Brunel University, England, 1994.

[8] A. M. Wazwaz, The numerical solution of fifth-order boundary-value
problems by the decomposition method, J. Comp. Appl. Math. 136 (2001)
259-270.

[9] D. J. Fyfe, Linear dependence relations connecting equal interval Nth
degree splines and their derivatives, J. Inst. Math. Appl. 7 (1971) 398-406.

[10] Siraj-ul-Islam, M.A. Khan, A numerical method based on non-polynomial
sextic spline functions for the solution of special fifth-order
boundary-value problems, Appl. Math. Comput. 181 (2006) 356-361.

[11] H. N. Caglar, S.H. Caglar, E.E. Twizell, The numerical solution of
fifth-order boundary-value problems with sixth degree B-spline functions,
Appl. Math. Lett. 12 (1999) 25-30.

[12] M. A. Noor, S.T. Mohyud-Din, Variational iteration technique for
solving fifth-order boundary-value problems, Appl. Math. Comput. 189 (2007)
1929-1942.

[13] M. A. Noor, S.T. Mohyud-Din, An efficient algorithm for solving
fifth-order boundary-value problems, Math. Comput. Modelling 45 (2007)
954-964.

[14] M. A. Khan, Siraj-ul-Islam, I.A. Tirmizi, E.H. Twizell, S. Ashraf, A
class of methods based on non-polynomial sextic spline functions for the
solution of a special fifth-order boundary-value problems, J. Math. Anal.
Applic. 321 (2006) 651-660.

[15] M. El-Gamel, Sinc and numerical solution of fifth-order boundary-value
problems, Appl. Math. Comput. 187 (2007) 1417-1433.

[16] A. Lamnii, H. Mraoui, D. Sbibih, A. Tijini, Sextic spline solution of
fifth-order boundary-value problems, Math. Comput. Simulation 77 (2008)
237-246.

[17] S. Siddiqi, G. Akram, Solution of fifth-order boundary-value problems
using nonpolynomial spline technique, Appl. Math. Comput. 175 (2006)
1574-1581.

[18] S. Siddiqi, G. Akram, Sextic spline solutions of fifth-order
boundary-value problems, Appl. Math. Lett. 20 (2007) 591-597.

[19] C. Wang, Z.X. Lee, Y. Kuo, Application of residual correction method in
calculating upper and lower approximate solutions fifth-order boundary-value
problems, Appl. Math. Comput. 199 (2008) 677-690.

[20] M. Cui, Z. Deng, Solutions to the definite solution problem of
differential equations in space $W_{2}^{l}\left[ 0,1\right] $, Adv. Math. 17
(1986) 327-328.

[21] M. Cui, Y. Lin, Nonlinear Numerical Analysis in the Reproducing Kernel
Spaces, Nova Science Publishers, New York, 2009.

[22] F. Geng, M. Cui, Solving a nonlinear system of second order boundary
value problems, J. Math. Anal. Appl. 327 (2007) 1167-1181.

[23] H. Yao, M. Cui, a new algorithm for a class of singular boundary value
problems, Appl. Math. Comput. 186 (2007) 1183-1191.

[24] W. Wang, M. Cui, B. Han, A new method for solving a class of singular
two-point boundary value problems, Appl. Math. Comput. 206 (2008) 721-727.

[25] Y. Zhou, Y. Lin, M. Cui, An efficient computational method for second
order boundary value problemsof nonlinear differential equations, Appl.
Math. Comput. 194 (2007) 357-365.

[26] X. L\"{u}, M. Cui, Analytic solutions to a class of nonlinear
infinite-delay-differential equations, J. Math. Anal. Appl. 343 (2008)
724-732.

[27] Y. L. Wang, L. Chao, Using reproducing kernel for solving a class of
partial differential equation with variable-coefficients, Appl. Math. Mech.
29 (2008) 129-137.

[28] F. Li, M. Cui, A best approximation for the solution of one-dimensional
variable-coefficient Burgers' equation, Num. Meth. Partial Dif. Eq. 25
(2009) 1353-1365.

[29] S. Zhou, M. Cui, Approximate solution for a variable-coefficient
semilinear heat equation with nonlocal boundary conditions, Int. J. Comp.
Math. 86 (2009) 2248-2258.

[30] F. Geng, M. Cui, New method based on the HPM and RKHSM for solving
forced Duffing equations with integral boundary conditions, J. Comp. Appl.
Math. 233 (2009) 165-172.

[31] J. Du, M. Cui, Solving the forced Duffing equations with integral
boundary conditions in the reproducing kernel space, Int. J. Comp. Math. 87
(2010) 2088-2100.

[32] X. Lv, M. Cui, An efficient computational method for linear fifth-order
two-point boundary value problems, J. Comp. Appl. Math. 234 (2010) 1551-1558.

[33] W. Jiang, M. Cui, Constructive proof for existence of nonlinear
two-point boundary value problems, Appl. Math. Comput. 215 (2009) 1937-1948.

[34] J. Du, M. Cui, Constructive proof of existence for a class of
fourth-order nonlinear BVPs, Comp. Math. Applic. 59 (2010) 903-911.

[35] M. Cui, H. Du, Representation of exact solution for the nonlinear
Volterra-Fredholm integral equations, Appl. Math. Comput. 182 (2006)
1795-1802.

[36] B.Y. Wu, X.Y. Li, Iterative reproducing kernel method for nonlinear
oscillator with discontinuity, Appl. Math. Lett. 23 (2010) 1301-1304.

[37] W. Jiang, Y. Lin, Representation of exact solution for the
time-fractional telegraph equation in the reproducing kernel space, Commun.
Nonlinear Sci. Numer. Simulat. 16 (2011) 3639-3645.

[38] F. Geng, M. Cui, A reproducing kernel method for solving nonlocal
fractional boundary value problems, Appl. Math. Comput. in press.

[39] Y. Lin, M. Cui, A numerical solution to nonlinear multi-point
boundary-value problems in the reproducing kernel space, Math. Meth. Appl.
Sci. 34 (2011) 44-47.

[40] F. Geng, A numerical algorithm for nonlinear multi-point boundary value
problems, J. Comp. Appl. Math. in press.

[41] M. Mohammadi, R. Mokhtari, Solving the generalized regularized long
wave equation on the basis of a reproducing kernel space, J. Comput. Appl.
Math. 235 (2011) 4003-4014.

[42] B. Y. Wu, X.Y. Li, A new algorithm for a class of linear nonlocal
boundary value problems based on the reproducing kernel method, Appl. Math.
Lett. 24 (2011) 156-159.

[43] H. Yao, Y. Lin, New algorithm for solving a nonlinear hyperbolic
telegraph equation with an integral condition, Int. J. Numer. Meth.
Biomedical Eng. 27 (2011) 1558-1568.

[44] J. Zhang, The numerical solution of fifth-order boundary value problems
by the variational iteration method, Comp. Math. Applic. 58 (2009) 2347-2350.

\end{document}